\documentstyle[12pt]{article}  
\begin{document}
\input amssym.def
\input amssym
\baselineskip15pt
\textwidth=12truecm
\textheight=20truecm
\hoffset-.1cm
\voffset-.5cm
\def\Hom{\hbox{\rm Hom}}
\def\End{\hbox{\rm End}}
\def\J{{\cal J}}
\def\I{{\cal I}}
\def\N{{\cal N}}
\def\M{{\cal M}}
\def\O{{\cal O}}
\def\H{{\cal H}}
\def\Aut{\hbox{\rm Aut}}
\def\K{{\cal K}}
\def\T{{\cal T}}
\def\X{{\cal X}}
\def\F{{\cal F}}
\def\L{{\cal L}}
\def\pp{\psi}
\def\ff{\varphi}
\def\A{{\cal A}}
\def\B{{\cal B}}
\def\C{{\cal C}}
\def\Z{{\cal Z}}
\def\D{{\cal D}}
\def\ov{\overline}
\def\h#1{\hbox{\rm #1}}
\font\got=eufm10 scaled\magstephalf
\font\eightrm=cmr8

\title{The $C^*$--algebra of a Hilbert Bimodule}
\author{\\ \\ Sergio Doplicher\\
Dipartimento di Matematica,\\
 Universit\`a di Roma ``La Sapienza'',\\
I--00185 Roma, Italy\\ \\
Claudia Pinzari\\
Dipartimento di Matematica,\\
Universit\`a di Roma ``Tor Vergata'',\\
I--00133 Roma, Italy\\ \\
Rita Zuccante\\
Dipartimento di Matematica,\\
Universit\`a di Firenze,\\
I--51134 Firenze, Italy\\ }
\date{}
\maketitle

\vfill
Research 
supported by  MURST,  CNR--GNAFA and European Community
\pagebreak

\begin{abstract}
   We regard a right Hilbert $C^*$--module $X$ over a $C^*$--algebra $\A$
endowed with an isometric $^*$--homomorphism $\phi: \A\to\L_\A(X)$ 
as an object $X_\A$ of the $C^*$--category of right Hilbert $\A$--modules.
   Following [11], we associate to it a $C^*$--algebra $\O_{X_\A}$ containing 
$X$ as a ``Hilbert $\A$--bimodule in ${\O}_{X_\A}$''. If $X$ is full and finite 
projective ${\O}_{X_\A}$ is the $C^*$--algebra $C^*(X)\ ,$ the generalization of the Cuntz--Krieger
algebras introduced by Pimsner [27]. More generally, $C^*(X)$ is canonically 
embedded in ${\O}_{X_\A}$ as the $C^*$--subalgebra generated by $X\ .$ Conversely, if $X$ is full 
${\O}_{X_\A}$ is canonically embedded in $C^*(X)^{**}\ .$ 

   Moreover, regarding $X$ as an object $_\A X_\A$ of the $C^*$--category of 
Hilbert $\A$--bimodules, we associate to it a $C^*$--subalgebra $\O_{_\A X_\A}$ 
of $\O_{X_\A}$ commuting with $\A\ ,$ on which $X$ induces a canonical 
endomorphism $\rho\ .$ We discuss conditions under which $\A$ and $\O_{_\A X_\A}$ 
are the relative commutant of each other and $X$ is precisely the subspace
of intertwiners in $\O_{X_\A}$ between the identity and $\rho$ on 
$\O_{_\A X_\A}\ .$

   We also discuss conditions which imply the simplicity of $C^*(X)$ 
or of $\O_{X_\A}\ ;$ in particular, if $X$ is finite projective
and full, $C^*(X)$ will be simple if $\A$ is $X$--simple and the 
``Connes spectrum'' of $X$ is ${\Bbb T}$.   
\end{abstract}

\section{Introduction}

Let $\C\subset\B$ be an inclusion of $C^*$--algebras and denote 
by $\A=\C'\cap\B$ the relative commutant. If $\rho$ is an endomorphism 
of $\C\ ,$ the subset $X_\rho$ of $\B$ defined by 
$$X_\rho=\{\psi\in\B| \psi C=\rho(C)\psi,\ C\in\C\}\eqno(1.1)$$
is a {\it Hilbert} $\A$--{\it bimodule in} $\B$, in the sense that $X_\rho$
is a closed subspace, stable under left and right multiplication
by elements of $\A$, and equipped with an $\A$--valued
right $\A$--linear inner product given by
$$<\psi, \psi'>_\A=\psi^*\psi'\ ,\quad \psi, 
\psi'\in X_\rho$$
such that $\|<\psi, \psi>_\A\|={\|\psi\|_\B}^2\ .$
We say that $\rho$ is {\it inner} in $\B$ if $X_\rho$ is  finite projective
as a right $\A$--module and if its  left annihilator in $\B$ is zero.

This notion reduces to that of inner endomorphism when, e.g., $\C=\B$
has centre ${\Bbb C}I\ ;$ if $\C\neq\B$ but $\A={\Bbb C}I$, $X_\rho$ is a 
Hilbert space 
in $\B$ and $\rho$ is the restriction to $\C$ of an inner endomorphism
of $\B$ [9, 10, 11], i.e. $\rho(C)=\sum_{1}^{d}\psi_iC{\psi_i}^*$, with
$\{\psi_i, i=1,\dots,d\}$ an orthonormal basis of $X_\rho\ .$ 

The crossed product of a unital $C^*$--algebra $\C$ with trivial centre
by the outer action of a discrete group [13, 19, 25] or by the action of a compact 
group dual [10] has the characteristic property that the objects 
(automorphisms,
resp. endomorphisms of $\C$) become inner in the crossed product $\B$,
and that $\A'\cap\B={\Bbb C}I\ .$

These notions of crossed products might prove too narrow to provide 
a scheme for an abstract duality theory of quantum groups in the spirit
of [11], or  for the related problem of describing the superselection structure 
of low dimensional QFT by a symmetry principle [12, 15]. In the last case, indeed,
no--go theorems indicate that the relative commutant of the observable 
algebra in the field algebra might have to be nontrivial [23, 29].

It is therefore interesting to study more general crossed products $\B$
associated to the pairs $\{\C, \rho\}$ and conditions ensuring existence 
and uniqueness, in particular of the $C^*$--algebra $\A$ appearing as the
relative commutant $\C'\cap\B\ .$ 

As a preliminary step towards this problem, that we hope to treat
elsewhere, we consider in this paper the situation where $X$ is given as
a Hilbert $C^*$--bimodule with coefficients in $\A$ (i.e. $X$ is
a right Hilbert $\A$--module with a monomorphism of $\A$ into the 
$C^*$--algebra $\L(X)$ of the adjointable module maps, defining the 
left action [27]).

With $X^r\ ,$ $r=0,1,2,\dots$ the bimodule tensor powers of $X$
(where $X^0=\A$ by convention) we can consider the following 
$C^*$--categories: 

--the strict {\it tensor} $C^*$--category $\T_X$
with objects $X^r\ ,$ $r\in{\Bbb N}_0$, and with arrows the adjointable 
right $\A$--module maps commuting with the left action of $\A$; 

--the $C^*$--category ${\cal S}_X$ with the same objects and with arrows
all the adjointable right $\A$--module maps. This is a strict
{\it semitensor} $C^*$--category in the sense that
on arrows only the tensor product on the right with the identity arrows of the category itself
is defined (cf. Section 2).

A general construction associates functorially to each object
$\rho$ in a strict tensor $C^*$--category a $C^*$--algebra 
$\O_\rho$ [11]. It is easy to verify that this applies without 
substantial modifications to objects in a strict {\it semitensor\/} $C^*$--category.
We can thus associate to the bimodule $X$ viewed as an object of 
${\cal S}_X$ (to mean this we will write for  short $X_\A$) a $C^*$--algebra 
$\O_{X_\A}$, where $\A$ is embedded as a $C^*$--subalgebra and
$X$ is embedded as a Hilbert $\A$--bimodule in $\O_{X_\A}\ .$
The $C^*$--algebra $C^*(X)$ constructed by Pimsner [27] from the 
bimodule $X$, generalizing the Cuntz--Krieger algebras, can be identified 
with the $C^*$--subalgebra of $\O_{X_\A}$ generated by $X$, and will coincide
with $\O_{X_\A}$ if $X$ is full and finite projective (Section 3).

The $C^*$--algebra $\O_{_\A X_\A}$ associated with $_\A X_\A\ ,$ i.e.
with  $X$ viewed as an object
of the {\it tensor\/} category $\T_X\ ,$ is embedded in the relative commutant
${\A}'\cap\O_{X_\A}$ and coincides with it if further conditions are
fulfilled (Proposition 3.4). $X$ induces a canonical endomorphism 
on ${\A}'\cap\O_{X_\A}$ which 
acts 
on $\O_{_\A X_\A}$  tensoring the arrows in $(X^r, X^s)$
with the identity arrow of $(X, X)$ on the left. 
We give conditions which guarantee that $\A$ is {\it normal} in 
$\O_{X_\A}\ ,$ i.e. $\A=(\A'\cap\O_{X_\A})'\cap\O_{X_\A}\ ;$ in this
case $X$ identifies with the $\A$--bimodule 
in $\B=\O_{X_\A}$ which induces $\rho$ on $\C=\A'\cap\O_{X_\A}$ in the 
sense of eq. (1.1).

If $X$ is full, $C^*(X)$ is the universal $C^*$--algebra 
containing $\A$ and $X$ as an $\A$--bimodule and generated by $X$;
$\O_{X_\A}$ can be canonically identified with a $C^*$--subalgebra of ${C^*(X)}^{**}$
(Theorem 3.3).

While $\O_{X_\A}$ generalizes the Cuntz algebras $\O_n$, $n<\infty$
when $X$ is finite projective, if $X$ is not it rather generalizes
the $C^*$--algebra $\O_H$ discussed in [6]. 

If $X$ is finite projective and full and $\A$ has no closed two sided
proper ideal $J$ such that $X^*JX\subset J$, then $C^*(X)$ is {\it simple}
if the Connes spectrum of the dual action of ${\Bbb Z}$ on the crossed
product of $C^*(X)$ with the canonical action of ${\Bbb T}$ is full,
i.e. coincides with ${\Bbb T}$. 
If furthermore there is a tensor power 
$X^s$ of $X$ containing an isometry which commutes with $\A$, then
$\O_{X_\A}$ is also simple. These and slightly more general
results are discussed in Section 4 (cf. Theorem 4.7).

\section{Representations of Hilbert Bimodules in $C^*$--Algebras}

A  strict {\it semitensor\/} $C^*$--category is a $C^*$--category
$\T$ for which the set of objects is a unital semigroup, with identity $\iota\ ,$
and such that for any object $\tau\in\T$  there is a $^*$--functor (``right tensoring'' with
the identity $1_\tau$ of $(\tau, \tau)$)
$$\Phi_\tau: (\rho, \sigma)\to(\rho\tau, \sigma\tau)\eqno(2.1)$$ such that
$$\Phi_\iota=id\ ,$$
$${\Phi_\omega}\circ{\Phi_\tau}=
{\Phi_{\tau\omega}}\ .$$
The product on the set of objects
will be referred to as the tensor product.
In other words $\Phi : \tau\to\Phi_\tau$ is a unital antihomomorphism
from the semigroup of objects of $\T$ to the semigroup $\End(\T)$
of $^*$--endofunctors of $\T\ .$ 
We will consider only cases where $\Phi_\tau$ is injective, and hence isometric.
Any strict tensor $C^*$--category is obviously semitensor 
choosing $\Phi_\tau: T\to T\times 1_\tau\ .$

Let $\A$ and $\B$ be $C^*$--algebras. A  Hilbert $\A$--$\B$--bimodule is a right
Hilbert $\B$--module $X$ 
(with $\B$--valued inner product denoted by $<x, y>_{\B}$)
 endowed with a faithful $^*$--homomorphism
$\phi: \A\to\L_\B(X)\ .$

It was shown in \cite{B} that a refinement of an argument by Dixmier on 
approximate units shows that if $X$ is countably generated as a 
 right Hilbert module then there exist elements $x_1, x_2, \dots$ of 
$X$ such that $\sum_{j=1}^{N}\vartheta_{x_j, x_j}$ is an approximate
unit for $\K_{\B}(X)\ ,$ the $C^*$--algebra of compact operators on $X\ .$
In particular every $x\in X$ is the norm limit
$$x=\sum_j\vartheta_{x_j, x_j}(x)=\sum_j x_j<x_j, x>_\B\ .$$ The set $\{x_j\}$ will be called a {\it basis} of 
$X\ .$  The use of a basis will be helpful to simplify our formalism, hence
throughout this paper
we will only consider countably generated Hilbert bimodules. However, 
most of our results extend to the more general setting. 

Let  $\B$ be a $C^*$--subalgebra of a $C^*$--algebra $\M\ .$ 
A right Hilbert $\B$--module contained in  $\M$ is a norm closed subspace
such that 
$$X\B\subseteq X\ ,\quad X^*X\subseteq\B$$
(for any pair of subspaces $X, Y\subset\M\ ,$ $XY$ denotes
  the closed linear subspace generated by operator products $xy\ ,$ $x\in X\ ,$ $y\in Y$).
If furthermore $\A\subset\M$ is a $C^*$--subalgebra satisfying
$$\A X\subseteq X\ ,\quad ax=0\ , x\in X\Rightarrow a=0\ ,$$
$X$ will be called a  Hilbert $\A$--$\B$--bimodule contained in $\M\ .$

If $X$ and $Y$ are respectively a  right Hibert $\B$--module and a Hilbert $\B$--$\C$--bimodule
in $\M$ then $XY$ is a right Hilbert $\C$--module in $\M$ naturally isomorphic to
 $X\otimes_{\B}Y\ .$ 

If $X$ and $Y$ are right Hilbert $\B$--modules in $\M$ then $YX^*$ is a subspace of
$\M$ naturally isomorphic
to the space $\K_\B(X, Y)$ of compact operators from $X$ to $Y\ .$ In general this identification
does not extend to the space
$\L_\B(X, Y)$ of $\B$--linear adjointable maps. However, $\L_\B(X, Y)$ may be 
recovered
as a subspace of $\M^{**}\ ,$ the enveloping von Neumann algebra of $\M\ .$
Let $1_X\in\M^{**}$ denote the identity of ${\overline{XX^*}^{uw}}\ ,$ the closure of
$XX^*$ in $\M^{**}$ in the ultraweak  topology.\medskip

\noindent{\bf 2.1. Proposition.} {\sl Let $X$ and $Y$ be right Hilbert $\B$--modules in $\M\ .$
Then setting 
$$(X, Y)_\B:=\{T\in {\M}^{**}: T1_X=1_YT=T\ , TX\subseteq Y, Y^*T\subseteq X^*\}$$
one defines a subspace of $\M^{**}\ ,$  in fact  contained in ${\overline{YX^*}^{uw}}\ ,$
which 
identifies naturally with $\L_B(X, Y)\ .$
If $X$ and $Y$ are Hilbert $\A$--$\B$--bimodules in $\M$
 then  $_\A(X, Y)_\B:=\A'\cap(X, Y)_\B$  corresponds in the above identification to 
the set of elements of $\L_\B(X, Y)$ that commute with the left 
$\A$--action.
}\medskip

\noindent{\bf Proof.} Any  $T\in(X, Y)_\B$ defines by multiplication in $\M^{**}$ an
 operator $\hat{T}: X\to Y$ with adjoint
$\hat{T^*}\ ,$ hence $\hat{T}\in\L_\B(X, Y)\ .$
Since $TXX^*\subseteq YX^{*}$ we conclude, approximating $1_X$ ultra strongly with elements of 
$XX^*\ ,$ that $T\in{\overline{YX^*}^{us}}\ .$ Furthermore $TX=0$ implies $T=0\ ,$ and this 
shows that $T\to\hat{T}$ is injective. On the other hand this map
is clearly isometric  from $YX^*$ to $\K_\B(X, Y)\ .$ If now $S\in\L_\B(X, Y)$  then for any basis 
$\{x_j\ , j=1, 2,\dots\}$ of $X\ ,$ $S\sum_{j=1}^{N}\vartheta_{x_j, x_j}$ is a norm bounded 
sequence of compact operators
hence it is of the form $\hat{T}_N\ ,$ with $T_N\in XX^*$ norm bounded and
 strictly convergent. Let $T\in{\overline{YX^*}^{uw}}$ be a weak limit point. Clearly $T1_X=1_YT=T\ .$
Furthermore for all $x\in X$
$T_Nx$ is norm convergent, necessarily to $Tx\ ,$ so $TX\subseteq Y\ .$  We also conclude
that $S=\hat{T}\ ,$ hence the map $T\to\hat{T}$ is surjective and the proof is complete.
$\hfill\square$\medskip

A representation of a $C^*$--category $\T$ in some $\B(H)$
is a collection of maps   $\F_{\rho, \sigma}: (\rho, \sigma)\to\B(H)\ ,$
 $\rho, \sigma\in\T$
such that for any pair of arrows $T\in(\rho, \sigma)\ ,$ $S\in(\sigma, \tau)\ ,$
$$\F_{\rho, \sigma}(T)^*=\F_{\sigma, \rho}(T^*)\ ,$$
$$\F_{\rho, \tau}(ST)=\F_{\sigma, \tau}(S)\F_{\rho, \sigma}(T)\ .$$
Let $\H_\B$ be the $C^*$--category of 
right Hilbert $\B$--bimodules: If $X$ and $Y$ are objects of $\H_\B$
the set of arrows from $X$ to $Y$  is $\L_\B(X, Y)$. Let $\T\subseteq\H_\B$
be a full subcategory.     Then the previous Proposition shows that if the objects of $\T$
embed in $\M$ as right Hilbert $\B$--modules then
there is a representation of $\T$ in the bounded linear operators on the 
Hilbert space of the universal representation of $\M\ .$

Note that in place of the universal representation  we may consider  
 any faithful representation  of $\M$ on some Hilbert space $H\ .$
 Indeed,   the subspace 
$(X, Y)_\B:=\{T\in\B(H): T1_X=1_YT=T\ ,TX\subseteq Y\ , Y^*T\subseteq X^*\}$
lies in $\overline{YX^*}^{uw}$ and again identifies naturally with $\L_\B(X, Y)$
($1_X$ is as before the identity of $\overline{XX^*}^{uw}\subseteq\B(H)$).
 It follows that there is still an obvious faithful representation
of $\T$ in $\B(H)\ .$ 

Our next aim is to extend the formalism of \cite{DRJOT} to Hilbert bimodules. We 
describe natural realizations of  categories of Hilbert bimodules faithfully represented in 
some $C^*$--algebra
 as endomorphism categories of a suitable $C^*$--algebra. Our starting point is 
the following. We   are given
 a unital semigroup  $\Delta$ of Hilbert bimodules over a $C^*$--algebra $\A$
contained in the $C^*$--algebra $\M\ .$ We  assume, for simplicity,
 that $\M$ is generated by the
elements of $\Delta\ .$ We  form the subspaces $(X, Y)_\A\ ,$ $X, Y\in\Delta\ ,$
 in $\M^{**}$ and  the category ${\cal S}_\Delta$ with arrows these intertwining spaces. We  denote
by $\tilde{\M}$ the $C^*$--subalgebra of $\M^{**}$ generated by the $(X, Y)_\A$'s.
It is now clear that ${\cal S}_\Delta$ is a strict semitensor $C^*$--category. If furthermore we define
$_\A(X, Y)_\A\subset(X, Y)_\A$ as the subspace of $\A$--bimodule maps, namely
$$_\A(X, Y)_\A=\{T\in(X, Y)_\A: aTx=Tax\ , a\in\A, x\in X\}\ ,$$ the  
 subcategory ${\cal T}_\Delta\subset{\cal S}_\Delta$ with the same objects of ${\cal S}_\Delta$ and arrows 
$_\A(X, Y)_\A\ ,$ is a strict tensor $C^*$--category.

Let $\B\subseteq\C$ be an inclusion of unital $C^*$--algebras, and let $\End_\C(\B)$ be the category
of  endomorphisms of $\B$ with arrows the intertwiners in $\C\ :$
$$(\rho, \sigma)=\{c\in\C: c\rho(I)=c\ , c\rho(b)=\sigma(b)c\ , b\in\B\}\ .$$
\medskip 

\noindent{\bf 2.2. Remark.} $\End_\C(\B)$ is a strict semitensor 
$C^*$--category by defining the tensor product on the set of objects
  to be the composition, and  
$\Phi_\tau: c\in(\rho, \sigma)\to c\in(\rho\tau, \sigma\tau)\ .$

 $\End_\B(\B)$ (simply denoted $\End(\B)$) is a tensor $C^*$--category by
$b\times b'=b\rho(b')\in(\rho\rho', \sigma\sigma')\ ,$ $b\in(\rho, \sigma)\ ,$ 
 $b'\in(\rho', \sigma')\ .$ 
\medskip

\noindent{\bf 2.3. Theorem.} {\sl Let $\Delta$ be a unital
semigroup of Hilbert $\A$--bimodules in a $C^*$--algebra $\M\ .$
With the above notation,  any $X\in\Delta$ induces a unique endomorphism
$\sigma_X$ on $\A'\cap\tilde{\M}$ such that 
$$\sigma_X(T)x=xT\ , x\in X\quad\ , T\in\A'\cap\tilde{\M}\ .$$
The map $X\in{\cal S}_\Delta\to\sigma_X\in\End_{\tilde{\M}}(\A'\cap\tilde{\M})$
that acts trivially on the arrows is a faithful functor of semitensor $C^*$--categories 
that restricts to a functor of tensor $C^*$--categories ${\cal T}_\Delta\to\End(\A'\cap\tilde{\M})\ .$ 
If furthermore $\A$ is normal in $\tilde{\M}$ then
 the images of these functors  are  full subcategories.}\medskip

\noindent{\bf Proof.} 
Let $\{x_1, x_2, \dots\}$ be a basis of $X\ .$ If $T\in{{\M}^{**}}^{+}$ then
the sequence of positive elements
 $\sum_{j=1}^{N}x_jT{x_j}^*$ is increasing and bounded in norm by $\|T\|\|1_X\|\ .$
Therefore $\sum_{j=1}^{N}x_jT{x_j}^*$ is strongly convergent to an element
$\phi(T)\in \M^{**}$ for any $T\in\M^{**}$ and $\phi$ is a norm $1$ positive map.
If $T\in(Y, Z)_\A\ ,$ for $Y, Z\in\Delta\ ,$ then clearly $\phi(T)\in \overline{XZY^*X^*}^{uw}$
and $\phi(T)XY\subseteq XZ$ and $\phi(T)^*XZ\subseteq XY\ ,$  hence 
$\phi(T)\in(XY, XZ)_\A\ .$ It follows that $\phi$ leaves
$\tilde{\M}$ globally invariant. Since $X^*X\subseteq\A\ ,$ $\phi$ is multiplicative
on $\A'\cap\tilde{\M}\ .$ Clearly if $T\in\A'\cap\tilde{\M}$ then 
$\sigma_X(T)x=xT$ for any $x\in X\ .$ Now   $\sigma_X(T)$ has support contained in  $1_X\ ,$ thus we
conclude that $\sigma_X(T)$ is independent on the basis. In particular, 
if $u$ is a unitary in 
${\A}$ (or in $\tilde{\A}:=\A+{\Bbb C}1_X$ if $\A$ does not have a unit) then the basis $\{ux_1, ux_2,
\dots\}$ induces the same map $\sigma_X\ ,$ thus $u$ commutes with $\sigma_X(\A'\cap\tilde{\M})\ ,$
i.e. $\sigma_X$ leaves $\A'\cap\tilde{\M}$ invariant. Finally, if $\A$ is normal in 
$\tilde{M}$ and $T\in(\sigma_X, \sigma_Y)$ then in particular  for any $x\in X$ and $y\in Y$
$y^*Tx\in{\A'\cap\tilde{M}}'\cap\tilde{M}=\A\ .$  For any basis $\{y_j\}$
of $Y\ ,$ $\sum_j y_j{y_j}^*Tx$ is norm converging to $Tx\ ,$ thus $TX\subset Y\ .$ Similarly, 
$T^*Y\subset X$, so $T\in(X, Y)_\A\ .$
$\hfill\square$\medskip

\section{The $C^*$--algebra $\O_\rho$} 
In this section we discuss the $C^*$--algebra $\O_\rho$ associated
with an object $\rho$ a strict semitensor $C^*$--category $\T\ .$ When we specialize
$\rho$ to a Hilbert bimodule $X\ ,$ $\A$ will be embedded in $\O_\rho$ as a 
subalgebra, and $X$ as a $\A$--bimodule. In view of Theorem 2.3 we will give 
sufficient conditions on $X$ in order that $\A$ is embedded as a normal subalgebra.

The construction of the $C^*$--algebras $\O_\rho$ was given in \cite{DRInv} 
 when $\rho$ is an object of a strict tensor $C^*$--category $\T\ .$
We are now interested, among others, in the categories ${\cal S}_X$ with objects the 
tensor powers of a bimodule $X$  and arrows  $(X^r, X^s)_\A\ ,$
$r, s\in{\Bbb N}_0\ ,$ so that ${\cal S}_X$ is only a strict semitensor $C^*$--category. 
However, the construction in \cite{DRInv} goes through without
substantial modifications and for the convenience of the reader  we sketch it here
in the case of a strict semitensor $C^*$--category. 

We first form the Banach space ${\O_{\rho}}^{(k)}$ inductive limit  of 
$(\rho^r, \rho^{r+k})$ via the maps 
$\Phi_\rho: (\rho^r, \rho^{r+k})\to (\rho^{r+1}, \rho^{r+k+1})\ .$ 
The composition 
and the $^*$--involution of $\T$ 
define on $\oplus_{k\in{\Bbb Z}}{\O_\rho}^{(k)}$ a structure 
of ${\Bbb Z}$--graded $^*$--algebra.  There is a unique $C^*$--norm  on 
$\oplus_{k\in{\Bbb Z}}{\O_\rho}^{(k)}$
 for which the automorphic
 action of ${\Bbb T}$ defined by the grading is isometric, and $\O_\rho$ is the completion
in that norm.
We denote by $^0\O_\rho$ the canonical dense $^*$--subalgebra
generated by images of intertwiners $(\rho^r, \rho^s)\ .$

If $\T$ is a genuine tensor $C^*$--category, 
tensoring
on the left by $1_{\rho}$ induces a canonical endomorphism, $\sigma_{\rho}$ of
$\O_\rho\ .$ 

Any $^*$--functor $\F: \T_1\to\T_2$ of strict semitensor $C^*$--categories
induces an obvious
 $^*$--homomorphism $\F_*: \O_\rho\to\O_{\F(\rho)}\ .$ 

Let $X$ be a Hilbert $\A$--bimodule as in Section 2. We will consider the semitensor 
$C^*$--category ${\cal S}_X$ with objects the $\A$--bimodule tensor powers $X^r$ of $X$ and arrows
the $(X^r, X^s)_\A\ ,$ the adjointable right $\A$--module  maps.
We will write $X_{\A}$ when $X$ is viewed as an object of this strict {\it semitensor\/}
 $C^*$--category. We can also consider the the strict {\it tensor\/} $C^*$--category $\T_X$
with the same objects and arrows the bimodule maps $_\A(X^r, X^s)_\A\ .$ We will
write $_\A X_{\A}$ when $X$ is considered as an object of this strict tensor category.

The construction of $\O_\rho$ applied to $\rho=X_{\A}$ yields a $C^*$--algebra 
$\O_{X_{\A}}$ that contains a copy of  $\A$ as embedded in $(X, X)_\A$ and
$X=\K_\A(\A, X)\subset(\A, X)_\A$ as a Hilbert $\A$--bimodule.
$\O_{X_{\A}}$ is generated as a Banach space by the $(X^r, X^s)_\A$'s and carries the action
$\alpha$ of ${\Bbb T}$ defined by the ${\Bbb Z}$--grading ${\O_{X_{\A}}}^{(k)}\ .$
\medskip

\noindent{\bf Remark.} The left annihilator of $X$ in $\O_{X_\A}$ is zero.
For, given $T\in{\O_{X_\A}}\ ,$ $Tx=0$ for all $x\in X$
implies, by Fourier analysis over the action $\alpha$
of ${\Bbb T}\ ,$ $T_kx=0$ for all $x\in X$,  $k\in{\Bbb Z}$, where $T_k$
is the projection of $T$ in ${\O_{X_\A}}^{(k)}\ .$ But each
${T_k}^*T_k$ can be approximated in norm by elements of $(X^r, X^r)_\A$ for
large $r$, and the norm on $(X^r, X^r)_\A$ is that of the corresponding
bounded operators on ${X^r}{X^r}^*\ .$ Thus $T_k=0$ and
$T=0\ .$
\medskip

\noindent{\bf Remark.} In the special case where $X$ is a 
Hilbert $\A$--bimodule in the $C^*$--algebra $\M\ ,$
$(X^r, X^s)_\A$ are identified as in Section 2 with the corresponding subspaces
of $\tilde{\M}\ ,$ but the closed linear span in $\tilde{\M}$ does not necessarily identify
with $\O_{X_{\A}}$  since the ${\Bbb Z}$--graded $^*$--subalgebra of $\tilde{\M}$
generated by the $(X^r, X^s)_\A$ does not necessarily carry
an automorphic action of $\T$ defined by the grading and continuous for the 
norm of $\tilde{\M}\ .$
\medskip

The following
result is an easy consequence of the definition of $\O_{_\A X_\A}$ and of
 functoriality  of the construction.\medskip

\noindent{\bf 3.1. Proposition.} {\sl Let $X$ and $Y$ be 
Hilbert $C^*$--bimodules over $C^*$--algebras $\A$ and $\B$ respectively,
and let $_\A\gamma_\B$ be a strong Morita 
equivalence such that $X$ and $\gamma Y\gamma^{-1}$  are isomorphic as Hilbert $C^*$--bimodules.
Then $\O_{_\A X_\A}$ and $\O_{_\B Y_\B}$ are isomorphic according to an isomorphism
that transforms $_\A(X^r, X^s)_\A$ into $_\B(Y^r, Y^s)_\B\ .$}\medskip 

Pimsner defined  in \cite{Pi} the universal $C^*$--algebra generated by a Hilbert 
bimodule $X\ .$ These $C^*$--algebras 
  are generalizations of the Cuntz--Krieger algebras and we shall refer to them as
CKP--algebras. In the following Proposition 
we relate the algebras  $\O_{X_{\A}}$ to the CKP--algebras.
\medskip

\noindent{\bf 3.2. Proposition.} {\sl Let $X$ be a  Hilbert $\A$--bimodule and
$C^*(X)$ the associated CKP--algebra. The identity map on $X$ extends to a $^*$--isomorphism
of $C^*(X)$ onto the $C^*$--subalgebra of $\O_{X_\A}$ generated by $X\ ,$
which is onto $\O_{X_\A}$ if $X$ is full and  projective.}\medskip

\noindent{\bf Proof.} Following Pimsner \cite{Pi}, we consider
 $\F(X)\ ,$  the full Fock space of $X\ ,$ and $J(\F(X))\ ,$ 
the $C^*$--subalgebra of $\L_\A(\F(X))$ generated by $\L_\A(\oplus_{n=0}^{p} X^n)\ ,$
$p\in{\Bbb N}\ .$ For any $x\in X\ ,$ let $S_x$ be the image in $M(J(\F(X)))/J(\F(X))$
of  the operator 
that tensors on the left by $x\ .$ The CKP--algebra 
 is  the $C^*$--subalgebra generated by $S_x\ , x\in X\ .$
The automorphic action $\beta$ of ${\Bbb T}$ on $\L(\F(X))$
implemented by  the unitary operators $U(z)$
on $\F(X)$ defined by $U(z)x=z^kx\ ,$ $x\in{X}^k\ ,$ $k\in{\Bbb N}_0\ ,$ induces
an action on the quotient $C^*$--algebra, that 
restricts to an action $\gamma$ on the CKP--algebra such that 
$\gamma_z(x)=zx$ for $x\in X\ .$
It follows that  the $^*$--subalgebra generated by $S_X$ is contained in $^0\O_{X_{\A}}$ in
a canonical way, and that this is an equality if $X$ is full and finite projective. Clearly, 
 the canonical action $\alpha$ corresponds to $\gamma\ .$ 
$\hfill\square$\medskip

\noindent{\bf 3.3. Theorem.} {\sl Let $X$ be a full Hilbert $\A$--bimodule contained
in $\M$
such that $\M$ is generated by $X$ as a $C^*$--algebra
(hence $\tilde{\M}$ is generated by the $(X^r, X^s)_\A$'s). 
The following are equivalent: \begin{description}
\item{ i)} $\M$ is the universal $C^*$--algebra with the properties above,
\item{ ii)} $\O_{X_\A}$ is canonically isomorphic to $\tilde{\M}\ ,$ i.e. there is a
$^*$--isomorphism acting as the identity on $(X^r, X^s)_\A\ ,$ $r, s\in{\Bbb N}_0\ ,$
\item{ iii)} the CKP--algebra $C^*(X)\subset\O_{X_\A}$ is canonically isomorphic to
$\M\ ,$ i.e. there is an isomorphism acting as the identity on $X\ ,$
\item{ iv)} there is an action $\alpha: {\Bbb T}\to\Aut(\M)$ such that $\alpha_z(x)=zx\ ,$
$z\in{\Bbb T}\ ,$ $x\in X\ .$
\end{description}
}\medskip

\noindent{\bf Proof.} If there is an action $\alpha$ as in {\sl iv)} then the  bitransposed action
${\alpha}^{**}: {\Bbb T}\to\Aut(\M^{**})$ restricts to 
an action on $\tilde{\M}$, still denoted by $\alpha\ ,$
   such that $\alpha_z(T)=z^{s-r}T\ ,$
$T\in(X^r, X^s)_\A\ ,$ and this shows the equivalence of {\sl ii)} with {\sl iv)} and with 
{\sl iii)} as well, 
in view of the previous Proposition. If {\sl i)} holds then {\sl iv)}
 follows from the universality property
of $\M\ .$ Finally, {\sl iii)}$\Rightarrow$ {\sl i)} was proved in \cite{Pi}.
$\hfill\square$\medskip

Theorem 3.3 can be easily reformulated without assuming that $X$ is full,
but  requiring that $\M$ is the $C^*$--algebra generated by $X$ and $\A\ .$
In this case, condition {\sl iii)} modifies requiring that
there is an isomorphism of $\M$ with the augmented algebra of Pimsner
\cite{Pi} which identifies the embeddings of $X\ ,$ respectively of $\A\ ,$
in those algebras. In condition {\sl iv)} the action $\alpha$ will be 
further required to be trivial on $\A\ .$

In view of condition {\sl i)} the CKP--algebra $C^*(X)\subset\O_{X_\A}$ can be thought of as
the crossed product of $\A$ by $X$ in the spirit of \cite{AEiEx} where, however, 
only bimodules of a more restricted class
were considered. \medskip

\noindent{\bf 3.4. Proposition.} {\sl\begin{description}
\item {\rm a)}
The  inclusion functor $\iota: {\cal T}_X\subset{\cal S}_X$ induces an inclusion 
$^*$--monomorphism
$$\iota_*:\O_{_\A X_{\A}}\to\O_{X_{\A}}$$ such that 
$$\iota_*(\O_{_\A X_{\A}})\subseteq\A'\cap\O_{X_{\A}}\ .$$
We have that $\sigma_X\circ\iota_*=\iota_*\circ\sigma_X\ .$
\item {\rm b)} If for some 
$s\in{\Bbb N}\ ,$ 
$_{\A}(\A, X^s)_\A$ contains an isometry  then 
$$\iota_*(\O_{_\A X_{\A}})={\A}'\cap\O_{X_\A}\ .$$
\end{description}}\medskip

\noindent{\bf Proof.} Part a) follows from the fact that the dual action of ${\Bbb T}$ on $\O_{X_{\A}}$
that defines its ${\Bbb Z}$--grading transform $_\A(X^r, X^s)_\A$ according to the character
$s-r\in{\Bbb Z}\ ,$ hence their linear span coincides with $^0\O_{_\A X_{\A}}\ .$
The canonical norm of $\O_{X_{\A}}\ ,$
i.e. the one for which the ${\Bbb T}$--action is isometric,
 restricts to the canonical norm of $^0\O_{_\A X_{\A}}\ .$ Since
 $\iota_*(\O_{_\A X_{\A}})$ and ${\A}'\cap\O_{\iota(X)}$ are globally invariant under the 
action of ${\Bbb T}\ ,$ to prove b) it suffices to show that the corresponding ${\Bbb T}$--eigenspaces are equal.
Let $R\in _\A(\A, X^s)_\A$ be an isometry. Using $R^{p+p'}={\sigma_X}^{ps}(R^{p'})R^p\ ,$ one can easily show that for $T$ in 
some $(X^p\ , X^{p+k})_\A$ the sequence ${\sigma_X}^{r+k}(R^{p'})^*T{\sigma_X}^r(R^{p'})\ ,$
$p'\in{\Bbb N}\ ,$ is eventually equal to a constant element
of $(X^r\ , X^{r+k})_\A\ .$ Thus the formula
$$E_r(T)=\lim_p {\sigma_X}^{r+k}(R^p)^*T{\sigma_X}^r(R^p)$$ 
defines a norm one  projection $E_r$ from ${\O_{_\A X_{\A}}}^{(k)}\ ,$  
the closure of $^0{\O_{_\A X_{\A}}}^{(k)}$ in $\O_{_\A X_{\A}}\ ,$ onto $(X^r\ , X^{r+k})_\A$ 
that acts identically on $(X^r\ , X^{r+k})_\A$ and satisfies 
$E_r(aTa')=aE_r(T)a'\ ,$ $a\ ,$ $a'\in{\A}\ .$ It follows that the sequence 
 $E_r$ is pointwise convergent to the identity map, thus if
 $T\in{\A}'\cap{\O_{X_\A}}^{(k)}$ then $E_r(T)\in(X^r, 
X^{r+k})_{\A}$ and approximates  $T\ .$
$\hfill\square$\medskip

The functorial properties of the construction of $\O_{X_\A}$ imply
that to each unitary $U\in _\A(X, X)_\A$ we can associate a canonical
automorphism $\sigma_U$ of $\O_{X_\A}\ ,$ leaving
$\O_{_\A X_\A}$ globally stable, such that
$$\sigma_U(x)=Ux\ ,\quad x\in X\ .$$
We thus establish an isomorphism between ${\cal U}(_\A(X, X)_\A)$ and the
group of all automorphisms of $\O_{X_\A}$ leaving $\A$ pointwise fixed and
$X$ globally stable.

The restriction to $\O_{_\A X_\A}$ of such an automorphism commutes with
$\sigma_X\ ;$ hence for each subgroup $G$ of ${\cal U}(_\A(X, X)_\A)$
the fixed point subalgebra ${\O_{_\A X_\A}}^G$ is globally stable under
$\sigma_X\ .$ Thus $\sigma_X$ induces an endomorphism $\sigma_G$ of 
${\O_{_\A X_\A}}^G\ .$

The systems (${\O_{_\A X_\A}}^G$, $\sigma_G$) have been extensively studied
when $\A={\Bbb C}\ ;$ we hope to turn to the general case where
$\A\neq{\Bbb C}$ and $G$ is replaced by a quantum group.

In the remaining part of this section we focus our attention on how $\A$ is embedded
in $\O_{X_\A}\ ,$ more precisely,
in view of Theorem 2.3 we look for conditions that $X$ should satisfy so that
$\A$ is normal in $\O_{X_{\A}}\ .$

A Hilbert $\A$--bimodule $X$ with left $\A$--action $\phi: \A\to\L_\B(X)$ is called  
nonsingular if
 $\vartheta_{x, x}\in{\phi(\A)}$ for some $x\in X$ implies $x=0\ .$
The trivial bimodule $\A$ is always singular.
It is easy to see that if $X$ is nonsingular then $Y\otimes_\A X$ is nonsingular for any Hilbert
$\A$--bimodule $Y\ .$
In particular, powers of nonsingular bimodules are  
nonsingular.

Let ${\A}$ be a unital,  purely infinite $C^*$--algebra, and let $X$ be 
a Hilbert $\A$--bimodule such that the left $\A$--action $\phi: \A\to\L_\B(X)$
 is unital.
 Then
$X$ is singular if and only if it is singly generated. In fact, if there is 
$a\in\A$ such that $\vartheta_{x, x}=\phi(a)\neq0\ ,$ then by the pure infiniteness
of $\A$ there is $b\in\A$ such that $b^*ab=I\ ,$ hence $\vartheta_{bx, bx}=I\ ,$
and this is to say that $X$ is singly generated. \medskip

\noindent{\bf 3.5. Proposition.} {\sl Let $X$ be a Hilbert $\A$--bimodule
in $\M\ .$ \begin{description}
\item{\rm a)}  If $X$ is nonsingular and
 $_\A(\A, X^s)_\A$ contains an isometry $S$ for some $s>1$ then 
$C^*(S)'\cap\O_{X_{\A}}=\A\ .$
 \item{\rm b)} If  there are isometries 
$S_k\in _\A(\A, X^{n(k)})_\A$ such that
$${S_k}^*{\sigma_X}^k({S_k})=\lambda_k$$
with $\|\lambda_k\|<1\ ,$ then $C^*(S_{k}\ , k=1, 2,\dots)'\cap\O_{X_{\A}}=\A\ .$
\end{description} 
In both cases $\A$ is normal in $\O_{X_{\A}}\ .$}
\medskip

\noindent{\bf Proof.} Let $\B$ denote one of the relative commutants described in
a) or in b), and $S\in _\A(\A, X^s)_\A$ an isometry.
 We  show that $\B\cap{\O_{X_{\A}}}^{(k)}$ is zero for $k\neq 0$
and that it is contained in ${\A}$ for $k=0\ .$
Let  $\Phi$ be a weak limit point of the sequence $\Phi_p(T)=(S^p)^*TS^p$ in some faithful 
representation of 
$\O_{X_{\A}}$ on a Hilbert space. Clearly $\Phi({\O_{X_{\A}}}^{(k)})\subseteq X^k\ ,$ $k\in{\Bbb N}\ ,$ and 
$\Phi(T)=T\ ,$ $T\in\B\ .$
 Hence $\B\cap{\O_{X_{\A}}}^{(k)}$ 
is contained in $X^k$ for $k\in{\Bbb N}_0\ .$ Let $T$ be an element of this subspace with 
$k>0\ .$ If $X$ is nonsingular then 
$TT^*\in\B\cap{\O_{X_{\A}}}^{(0)}\subset\A\ ,$ so  $T=0\ ,$ and a) holds. To prove b), we note that  
$$T={S_k}^*TS_k={S_k}^*{\sigma_X}^k(S_k)T=\lambda_k T\ ,$$
thus $T=0\ .$
$\hfill\square$\medskip

As a  consequence of the previous result we can show normalcy of 
$\A$ in $\O_{X_{\A}}$ when $X$ is a real or pseudoreal bimodule with 
  dimension $>1$ in the sense of  \cite{LR}. More explicitly, and slightly
more generally, we have the following result.\medskip

\noindent{\bf 3.6. Corollary.} {\sl If there is an isometry $S\in _\A(\A, X^2)_\A$
such that $$\|S^*\sigma_X(S)\|<1$$ then $C^*({\sigma_X}^k(S), k=0,1,2,\dots)'\cap\O_{X_{\A}}=\A\ ,$ hence
$\A$ is normal in $\O_{X_{\A}}\ .$}\medskip

\noindent{\bf Proof.}
The isometries
$$S_k:={\sigma_X}^{k-1}(S)\dots\sigma_X(S)S\in _\A(\A, X^{2k})_\A$$ 
satisfy
${S_k}^*{\sigma_X}^k(S_k)=\lambda^k\ ,$ with $\lambda=S^*\sigma_X(S)\ .$
$\hfill\square$\medskip

\bigskip

\section{The Ideal Structure of $\O_{X_{\A}}$}
In the first part of this section we   introduce a natural class of $^*$--re\-pre\-senta\-tions
$\pi: \O_{X_{\A}}\to\B(H)\ ,$ called locally strictly continuous, and we generalize Pimsner's universality
result to the algebras $\O_{X_{\A}}\ .$ We then associate to $X$
certain Connes spectra which allow us to characterize simplicity of
${\O_{X_\A}}$ and of the CKP--algebra $C^*(X)\subset\O_{X_\A}$
in terms of a   suitable class of ideals
of $\A\ .$

The following  is a variant of Pimsner's universality result
to the $C^*$--algebras $\O_{X_{\A}}\ .$\medskip

\noindent{\bf 4.1. Theorem.} {\sl Let  $Y$ be a Hilbert bimodule over a 
$C^*$--algebra $\B$
  in $\B(H)\ ,$ and ${\cal D}$ the $C^*$--subalgebra of $\B(H)$ generated
by the subspaces $(Y^r, Y^s)_\B\ ,$ $r, s\geq0$. Assume that the left annihilator
of $Y$ in $\D$ is zero, and let
   $(U, \phi)$ be a pair consisting of   $^*$--isomorphism
$\phi: \A\to\B$ and a linear surjective 
  map $U: X\to Y$  which satisfies 
$$U(x)^*U(x')=\phi(x^*x')\ ,$$
$$U(xa)=U(x)\phi(a)\ ,\quad U(ax)=\phi(a)U(x)\ ,$$
for $a, a'\in\A\ ,$ $x, x'\in X\ .$ 
Then there is a  unique $^*$--representation $\pi: C^*(X)\to\B(H)$ 
that maps  $x\in X$ to $U(x)$, as in [27, Theorem 3.12], which furthermore
extends 
 to a unique $^*$--representation $\tilde{\pi}: \O_{X_\A}\to\B(H)$
via
$$\tilde{\pi}(T)\pi(A)=\pi(TA)\ ,\quad A\in X^s{X^r}^*\ ,\quad T\in(X^s, X^t)_\A\ ,$$
$$\tilde{\pi}(a)=\phi(a)\ , \quad a\in\A\ .$$
If ker$\pi$ is ${\Bbb T}$--invariant
then $\tilde{\pi}$ is faithful.}\medskip

\noindent{\bf Proof.} It is easy to see that for any
$T\in(X, X)_\A$ there is a unique operator $\pi_U(T)\in(Y, Y)_\B$
such that $\pi_U(T)Ux=U(Tx)$, $x\in X$, and that 
 $\pi_U$ is a $^*$--homomorphism s.t. $\pi_U(xy^*)=U(x)U(y)^*\ ,$ $x, y\in X$.
 Let $\{x_1, x_2,\dots\}$ be a basis of $X$. Since $U$
has dense range,  $\{U(x_1), U(x_2), \dots\}$ is a basis of $Y$.
Since the left annihilator of $Y$ in the $C^*$--subalgebra
$C^*(Y, \B)$ of $\B(H)$ generated by $Y$ and $\B$ is zero, for any $a\in\A\cap XX^*$, 
 $\sum_i\phi(a){U(x_i)U(x_i)^*}=\sum_i\pi_U(ax_i{x_i}^*)$
 is norm converging to $\phi(a)$, therefore 
by [27, Theorem 3.12]
there is a unique $^*$--representation $\pi$ of
 $C^*(X)\subseteq\O_{X_{\A}}$ on $H$ such that $\pi(x)=U(x)\ .$ 
Now if max$\{r,s\}>0\ ,$
 the restriction of $\pi$ to  $X^s{X^r}^*$
 extends uniquely to a map $\tilde{\pi}_{r, s}: (X^r, X^s)_\A\to
(Y^r, Y^s)_\B\subset\B(H)$ 
such that for
$T\in(X^s, X^t)_\A$, $A\in{X^s}{X^r}^*$, $B\in X^v{X^t}^*$,
  $$\tilde{\pi}_{s, t}(T)\pi(A)=\pi(TA)\ ,$$
$$\pi(B)\tilde{\pi}_{s, t}(T)=\pi(BT)\ .$$ 
We set, by convention, $\tilde{\pi}_{0, 0}=\phi:\A\to\B\subset\B(H)\ .$
Uniqueness implies $\tilde{\pi}_{s, t}(T)^*=\tilde{\pi}_{t, s}(T^*)\ ,$
 $\tilde{\pi}_{s, t}\tilde{\pi}_{r, s}=\tilde{\pi}_{r, t}$, 
and also that the restriction of $\tilde{\pi}_{s+1, t+1}$ to
$(X^s, X^t)_\A$ coincides with $\tilde{\pi}_{s, t}$ since 
the left annihilator of $Y$ in $C^*(\{(Y^r, Y^s)_\B, r,s\geq0)\subset\B(H)$ is zero. 
We can thus define a unique   $^*$--homomorphism 
$\tilde{\pi}: {^0\O}_{X_\A}\to\B(H)$ extending $\tilde{\pi}_{r, s}$ on $(X^r, X^s)_\A$.
 We show that $\tilde{\pi}$ is norm continuous.
Let  $T=\sum_{k=-n}^{n} T_k$ be an element of $^0\O_{X_{\A}}\ ,$
with $T_k\in(X^r\ , X^{r+k})_\A$ for a suitable $r$
and $k=-n,\dots, n$ and let $1_{F}$ be the support of a finitely
generated right $\A$--submodule of $X^r\ ,$ so $T1_{F}\in C^*(X)\ .$ Then
$$\|\pi(T1_{F})\|\leq\|T1_F\|\leq\|T\|$$ for all $F$ implies 
$\|\tilde{\pi}(T)\|\leq\|T\|\ .$

Assume now that ker$\pi$ is globally invariant under the action of
${\Bbb T}\ .$ 
Then $\tilde{\pi}_{r, r}$ is faithful on $(X^r, X^r)_\A$ since the left annihilator
of $X^r$ in $\O_{X_\A}$ is zero, therefore,
since  ker$\pi\cap{\O_{X_{\A}}}^0$ is the inductive limit
of ker$\pi\cap(X^r, X^r)_\A\ ,$
 $\tilde{\pi}$ is faithful
 on ${\O_{X_{\A}}}^0$, hence, being ker$\tilde{\pi}$ ${\Bbb T}$--invariant,
 $\tilde{\pi}$ is faithful on $\O_{X_\A}\ .$
$\hfill\square$\medskip

As a consequence of Theorem 4.1 the correspondence between unitaries
and endomorphisms of the Cuntz algebras generalizes as follows.

\medskip

\noindent{\bf 4.2. Proposition.} {\sl
 Any unitary $U\in{\A}'\cap{\O_{X_{\A}}}$ defines an endomorphism $\lambda_U$
of $\O_{X_{\A}}$ acting trivially on ${\A}$ by
$$\lambda_U(x)=Ux\ ,\quad x\in X\ .$$
If $U\in{\A}'\cap{\O_{X_{\A}}}^0$ then $\lambda_U$ is a monomorphism.

If $X$ is finite projective, the correspondence $U\to\lambda_U$
is a one to one map of the unitaries in $\A'\cap\O_{X_\A}$ onto the 
endomorphisms 
of $\O_{X_\A}$ leaving $\A$ pointwise fixed, which extends the canonical
action of ${\cal U}(_\A(X, X)_\A)$ (cf. Section 3).}\medskip

\noindent{\bf Proof.} We represent $\O_{X_\A}$ faithfully on a Hilbert space 
$H$.
We have already noted that the left annihilator 
of $X$ in $\O_{X_\A}$ is zero (cf. a remark in Section 3), therefore also the left
annihilator  of $Y:=UX$  in 
$\O_{X_\A}$ (regarded as 
a Hilbert $\A$--bimodule in $\O_{X_\A}$)  
is zero. By Theorem $4.1$ there is a  unique $^*$--representation 
$\lambda_U$ of $\O_{X_{\A}}$ on $H$ such that $\lambda_U(x)=Ux$, $x\in X$
and acting trivially on $\A$ provided we show that the left annihilator
of $Y$ in 
$C^*\{(Y^r, Y^s)_\A, r,s\geq0\}\subset\B(H)$ 
is zero.
Now
$$Y^s{Y^r}^*=U\sigma_X(U)\dots{\sigma_X}^{s-1}(U)X^s{X^r}^*
{\sigma_X}^{r-1}(U^*)\dots
\sigma_X(U^*)U^*\ ,$$ therefore $$(Y^r, Y^s)_\A=
U\sigma_X(U)\dots{\sigma_X}^{s-1}(U)(X^r, X^s)_\A{\sigma_X}^{r-1}(U^*)\dots
\sigma_X(U^*)U^*
\subset\O_{X_\A},$$ and the claim follows
from the previous remarks. Note also that
 if $T\in X^s{X^r}^*$,
$$\lambda_U(T)=U\sigma_X(U)\dots{\sigma_X}^{s-1}(U)T{\sigma_X}^{r-1}(U^*)\dots
\sigma_X(U^*)U^*\ ,$$
therefore 
the same formula must hold for $T\in(X^r, X^s)_\A$, and we
conclude that $\lambda_U$
is an endomorphism of $\O_{X_{\A}}\ .$  If $U$ is a ${\Bbb T}$--fixed point
then $\lambda_U$ commutes with $\alpha\ ,$ so ker$\lambda_U$
is ${\Bbb T}$--invariant. 

Since the left annihilator of $X$ in ${\O_{X_\A}}$ is zero 
the map $U\to\lambda_U$ is one to one. If $X$ is finite projective 
and $x_1,\dots, x_d$ is a basis in $X$, for each endomorphism
$\lambda$ leaving $\A$ pointwise fixed we can define, following Cuntz,
$$U:=\sum_i \lambda(x_i){x_i}^*\ ,$$
so that $U$ is unitary. For $a\in\A\ ,$ 
$x\in X\ ,$ we have $$Uax=\lambda(ax)=a\lambda(x)=aUx$$
so that $U\in\A'\cap\O_{X_\A}$ and $\lambda=\lambda_U\ .$ If $\lambda(X)=X$
clearly $U\in_\A(X, X)_\A\ .$

$\hfill\square$\medskip

Our next aim is to determine the ideal structure of
$\O_{X_{\A}}$ in certain  cases of interest for our pourposes.
We first look at ideals invariant under the canonical action of
the circle group. 
 Let $\J$ be a closed ideal of  $\O_{X_{\A}}\ .$
 We call $\J$   locally strictly closed
if whenever one of $r$ and $s$ is 
nonzero $\J_{r, s}:=\J\cap(X^r, X^s)_\A$ is strictly closed in $(X^r, X^s)_{\A}\ .$
Note that in
 this case, $\J_{r, s}$ is the strict closure of $X^s\J\cap\A{X^r}^*$ in 
$(X^r, X^s)_{\A}\ .$
An ideal $J$ of $\A$ is called $X$--invariant if $X^*J X\subset J\ .$
As in \cite{KWP} we  associate to $J$ the ideal $J_X:=\{a\in\A: X^*aX\subset J\}$ which is 
a closed $X$--invariant ideal containing $J\ .$ We call $J$ 
$X$--{\it saturated\/}
if $J_X=J\ .$ Note that the zero ideal is $X$--saturated, and that, if 
$X$ is full and nondegenerate (in the sense that $\A X=X$) and if $J$ is proper 
then $J_X$ is proper.

\medskip

\noindent{\bf 4.3. Lemma.} {\sl 
\begin{description}
\item{\rm a)} Any ${\Bbb T}$--invariant  closed ideal
$\J$ of ${\O_{X_{\A}}}$ is the closed linear span of
 $\J_{r, s}\ ,$ $ r, s=0, 1, 2, \dots\ .$ Therefore,
if $\J$ is also l.s.c, it is determined by $\J\cap\A\ .$

\item{\rm b)} Let $J$ be an $X$--invariant, $X$--saturated ideal
of $\A\ ,$ and  let 
   $\tilde{\J}$ denote  the c.l.s. in $\O_{X_\A}$ of the
strict closures of $X^sJ{X^r}^*$ in $(X^r, X^s)_\A\ .$ If $X$ is full,
let $\J$ be 
 the c.l.s. of the $X^sJ{X^r}^*\ .$ Then $\tilde{\J}$ and $\J$ are respectively
 a locally strictly closed ${\Bbb T}$--invariant
ideal of $\O_{X_\A}$ and a closed ${\Bbb T}$--invariant ideal
of $C^*(X)$ such that $\tilde{\J}\cap\A=\J\cap\A=J\ .$
\end{description}}\medskip

\noindent{\bf Proof.}
We first note that if $\B$ is any $C^*$--algebra 
andowed with a continuous automorphic action $\alpha$ of ${\Bbb T}$
 and $\I$ and $\J$ are  closed
$\alpha$--invariant ideals of $\B$ such that the fixed point subalgebras
coincide: $\I^{\alpha}=\J^{\alpha}$ 
then $\I=\J\ .$ Indeed, by Fourier analysis $\I$ is generated by
by the subspaces $\I^{(k)}$ that transform like the character $k\in{\Bbb Z}=\hat{{\Bbb T}}\ .$
Furthermore by [26, Proposition 1.4.5] any element $T\in\I^{(k)}$ can be written in the 
form $T=u(T^*T)^{1/4}$ with $u\in\I\ ,$ hence $T\in\I \J^{\alpha}\subset\J\ ,$ i.e.
$\I\subset\J\ .$
Exchanging the role of $\I$ and $\J$ we deduce that $\J=\I\ .$
Let now $\J$ be a closed ${\Bbb T}$--invariant ideal of $\O_{X_\A}$
and let $\I$ be the closed linear span of $\J_{r, s}\ ,$ which is still
a ${\Bbb T}$--invariant ideal.
Since  the homogeneous part of 
$\O_{X_\A}$ is the inductive limit of 
$(X^r, X^r)_\A\ ,$   $\J$
is generated by the subspaces $\J_{r, r}$'s, hence $\I^{(0)}=\J^{(0)}\ ,$ therefore
the previous argument
shows that $\J$ is generated by the $\J_{r,s}\ .$
To prove b) we consider, for any $r\geq 0\ ,$ the ideal
$J_r$ of $(X^r, X^r)_\A\subset\O_{X_\A}$ 
defined by the strict closure of $X^rJ{X^r}^*$ in $(X^r, X^r)_\A\ ,$
 so that the inductive limit of the
$J_r$'s generates $\tilde{\J}\cap{\O_{X_\A}}^{(0)}\ .$ If $a\in\A\cap\tilde{\J}$ then clearly
$\lim_r dist(a, J_r\cap\A)=\lim_r dist(a, \A\cap\tilde{\J})=0\ .$ 
On the other hand $J_r\cap\A=J$
for all $r$ since $J$ is $X$--saturated, therefore $a\in J\ .$
It follows easily that  $\tilde{\J}\cap(X^r, X^r)_\A=J_r\ ,$ hence 
$\tilde{\J}$ is locally strictly closed and, clearly, ${\Bbb T}$--invariant.
In the second case, we may argue in the same way, replacing 
$\O_{X_\A}$ by $C^*(X)\ ,$ $\tilde{\J}$ by $\J\ ,$
$(X^r, X^r)_\A$ by $\A+XX^*+\dots X^r{X^r}^*\subset C^*(X)$ and  $\J_r$
by $J+XJX^*+\dots X^rJ{X^r}^*\ .$ Since $\A\cap J+XJX^*+\dots X^rJ{X^r}^*\subset
J_{X^r}=J\ ,$ we deduce as above that if $a\in\A\cap\J$ then $a\in J\ .$ 
$\hfill\square$\medskip

If $\J$ is a l.s.c. ideal of $\O_{X_\A}$ then
$\J\cap\A$ is always $X$--saturated. However, this is 
not necessarily true if $\J$ is an ideal of $C^*(X)\ .$
Indeed, this condition may be stated equivalently requiring that
if $\pi: C^*(X)\to C^*(X)/\J$ is the quotient map and $P$
is the support of the right $\pi(\A)$--module $\pi(X)$ contained in
$\pi(C^*(X))$ (hence $P\in\pi(C^*(X))^{**})$) then $\pi(a)P=0$ with $a\in\A$
implies $\pi(a)=0\ .$ In certain cases, e.g.  $\A\subset XX^*\ ,$
 then  $\J\cap\A$ is $X$--saturated for 
every closed ideal $\J$ of $C^*(X)\ .$ If 
 some positive power $X^s$ of $X$ contains an 
isometry commuting with $\A$ then every $X$--invariant ideal
is automatically $X$--saturated.
\medskip

\noindent{\bf 4.4. Proposition.} {\sl Let $J\to\J$ and $J\to\tilde{\J}$
be the maps described in the previous Lemma.\begin{description} 
\item{\rm a)}  $J\to\tilde{\J}$ 
is a  bijective correspondence between
$X$--invariant, $X$--saturated
ideals of ${\A}\ ,$  and 
${\Bbb T}$--invariant l.s.c. ideals of $\O_{X_{\A}}$ with
inverse $\tilde{\J}\to\tilde{\J}\cap\A\ .$

\item{\rm b)}
If $X$ is full and $\A\subset XX^*\ ,$ 
  $J\to\J$ is a bijective correspondence 
between the class of ideals of $\A$ described in
a) and the set of closed ${\Bbb T}$--invariant ideals of $C^*(X)$ 
with inverse the map $\J\to\J\cap\A\ .$\end{description}}\medskip

\noindent{\bf Proof.}
By  Lemma 4.3  and the above remarks we need only to show that
if $\A\subset XX^*$ then every closed ${\Bbb T}$--invariant 
ideal $\J$ of $C^*(X)$ is the c.l.s. of the subspaces
$X^r\J\cap\A{X^r}^*\ .$ Now $\A\subset XX^*$ implies 
that $X^r{X^r}^*\subset X^{r+1}{X^{r+1}}^*$ for all $r\in{\Bbb N}\ ,$
hence the homogeneous part of $C^*(X)$ is the inductive limit of 
$X^r{X^r}^*\ ,$ $r\in{\Bbb N}\ ,$ and this implies that
the homogeneous part of $\J$ is the inductive limit of 
$\J\cap X^r{X^r}^*=X^r\J\cap\A{X^r}^*\ ,$ therefore
 $\J$ is generated by the subspaces $X^s\J{X^r}^*\ .$
$\hfill\square$\medskip

We call $\A$ $X$--{\it simple\/} if it has no proper $X$--invariant, $X$--saturated
ideal, and $X$--{\it prime\/} if it has no pair of nonzero orthogonal $X$--invariant,
$X$--saturated ideals.\medskip

\noindent{\bf 4.5. Corollary.} {\sl If $X$ is a  Hilbert $\A$--bimodule,
the following properties are equivalent,
\begin{description}
\item{\rm a)} $\A$ is $X$--simple
(resp. $\A$ is $X$--prime),
\item{\rm b)}  $\O_{X_{\A}}$ has no proper locally strictly closed
${\Bbb T}$--invariant ideal (resp. $\O_{X_\A}$ has no pair of nonzero
orthogonal, locally strictly closed, ${\Bbb T}$--invariant ideals),\end{description}
Consider the following conditions:\begin{description}
\item{\rm i)} $\A\subset XX^*\ ,$
\item{\rm ii)} for some $s\in{\Bbb N}\ ,$ $X^s$ contains an isometry $S$
commuting with $\A\ .$\end{description}
If either i) or ii) holds and $X$ is full, a) and b) are also equivalent to
\begin{description}
\item{\rm c)} $C^*(X)$ is ${\Bbb T}$--simple (resp. $C^*(X)$ is ${\Bbb T}$--prime),
\end{description}
If ii) holds, a) and b) are equivalent to \begin{description}
\item{\rm d)} $\O_{X_\A}$ is ${\Bbb T}$--simple
(resp. $\O_{X_{\A}}$ is ${\Bbb T}$--prime).
\end{description}}\medskip

\noindent{\bf Proof.} We prove only the statements concerning simplicity,
those concerning primeness can be proved with similar arguments. The equivalence
of a) and b), and of a) and c), in the  case that i) holds, follow from Proposition
4.4. Note that by Lemma 4.3 c)$\Rightarrow$a) (even without assuming that
ii) holds). Conversely, assume that a) and ii) hold. Let $\J$ be
a nonzero ${\Bbb T}$--invariant ideal of $C^*(X)\ ,$ then $\J\cap\A$ is a nonzero,
 $X$--invariant,
 $X$--saturated ideal of $\A\ ,$ hence  $\J\cap\A=\A\ ,$ that implies $\J=C^*(X)\ .$
We are left to show that ii) and b) imply d).
Let $\J$ be a proper ${\Bbb T}$--invariant 
ideal of $\O_{X_{\A}}\ ,$  and define  
 $\tilde\J$ as the c.l.s.  of the strict closures of
$\J\cap(X^r, X^s)\ .$ $\tilde\J$ is a ${\Bbb T}$--invariant
ideal containing $\J\ .$ We claim  that $\tilde{\J}$ is locally strictly closed,
or, more precisely, that $\tilde\J\cap(X^r, X^s)$ is the strict closure
of ${\J\cap(X^r, X^s)}$ and that $\tilde\J\cap\A=\J\cap\A\ .$
 It suffices to prove the second assertion. Let $a$ be an element
of $\tilde\J\cap\A$ and let $T$ be in the strict closure of
some $\J\cap(X^{rs}, X^{rs})$ such that $\|a-T\|<\varepsilon\ , $ then
$\|a-{S^r}^*TS^r\|<\varepsilon\ ,$ hence $a\in\J\ .$ It follows, by b),
 that $\tilde\J=\O_{X_{\A}}\ .$
Let $T_{\alpha}$ be a net in some $\J\cap(X^{rs}, X^{rs})$ strictly converging 
to the identity, then ${S^r}^*T_{\alpha}S^r$ is a norm converging 
sequence in $\J$  to the identity, so $\J=\O_{X_{\A}}$ and the proof is complete.
$\hfill\square$\medskip

We denote by $\Gamma(X)$ and $\hat{\Gamma}(X)$ 
 the Connes spectra of the dual action 
$\hat\alpha$ of ${\Bbb Z}$ on $C^*(X)\times_{\alpha}{\Bbb T}$ and
$\O_{X_{\A}}\times_{\alpha}{\Bbb T}$ respectively.
By \cite{OP} (cf. also  Lemma 8.11.7 of \cite{Ped})  
$$\Gamma(X)=\{\lambda\in{\Bbb T}: \I\cap\alpha_{\lambda}(\I)\neq\{0\}\ 
, \I\  {\hbox{\rm  all closed non--zero ideal  of }} C^*(X)\}\ ,$$
$$\hat{\Gamma}(X)=\{\lambda\in{\Bbb T}: \I\cap\alpha_{\lambda}(\I)\neq\{0\}\ 
, \I\  {\hbox{\rm  all closed non--zero  ideal  of }} \O_{X_{\A}}\}\ .$$

We note that if $X$ is a Hilbert $\A$--bimodule such that $\O_{X_\A}\ ,$
(resp. $C^*(X)$)  
is prime or simple then clearly
$\hat{\Gamma}(X)={\Bbb T}$ (resp. $\Gamma(X)={\Bbb T}$). 
 Furthermore, by Lemma 4.3 $\A$ (resp. the $C^*$--subalgebra of $\A$
generated by the scalar products if $X$ is not full)
is necessarily $X$--prime or $X$--simple. 
The following 
results are a partial converse.\medskip 

\noindent{\bf 4.6. Proposition.} {\sl 
Let $X$ be  a  Hilbert $\A$--bimodule with $\A$  $X$--prime.
\begin{description}
\item{\rm a)} If $X$ is full and
one of the conditions i) or ii) of 4.5 is satisfied
and furthermore  and $\Gamma(X)={\Bbb T}$ then $C^*(X)$ is prime.
\item{\rm b)} If ii) of 4.5 is satisfied and $\hat{\Gamma}(X)={\Bbb T}$ then
$\O_{X_\A}$ is prime.
\end{description}}\medskip

\noindent{\bf Proof.} If $C^*(X)$ were not prime  then  the 
arguments that prove  
$(ii)\Rightarrow(i)$ of Theorem 8.11.10 in \cite{Ped} would
prove the existence of two 
non--zero 
${\Bbb T}$--invariant orthogonal  ideals in $C^*(X)\ ,$ but this is 
impossible because by 4.5 $C^*(X)$ is ${\Bbb 
T}$--prime.
We prove the second part of the Proposition. 
 $\hat{\Gamma}(X)={\Bbb T}$ and $\O_{X_{\A}}$ nonprime
 imply the existence of two orthogonal
${\Bbb T}$--invariant proper ideals of $\O_{X_{\A}}\ $ hence the existence of two 
proper orthogonal
$X$--invariant 
ideals of $\A$ again by 4.5.
$\hfill\square$\medskip

The above Proposition can be used to prove the following result.\medskip

\noindent{\bf 4.7. Theorem.} {\sl 
Let $X$ be  a  Hilbert $\A$--bimodule with $\A$  $X$--simple.
\begin{description}
\item{\rm a)} If $X$ is full and
one of the conditions i) or ii) of 4.5 is satisfied
and furthermore  and $\Gamma(X)={\Bbb T}$ then $C^*(X)$ is simple.
\item{\rm b)} If ii) of 4.5 is satisfied and $\hat{\Gamma}(X)={\Bbb T}$ then
$\O_{X_\A}$ is simple.
\end{description}}\medskip

\noindent{\bf Proof.} By Lemma 8.11.11 of 
\cite{Ped} it suffices
 to check that our assumptions in a) and b) imply primeness and ${\Bbb T}$--simplicity
of $C^*(X)$ and $\O_{X_\A}$ respectively, and this follows from Proposition
4.6 and Corollary 4.5.
  $\hfill\square$\medskip

 \noindent{\bf 
Acknowledgements} 
 C. P. 
takes pleasure to thank Gert K. Pedersen for warm hospitality and discussions
while she was visiting the Mathematics Institute of the University of Copenhagen
under the Human Capital and Mobility Program.
 \bigskip

\end{document}